\documentclass{article}

\usepackage{amssymb}
\usepackage{amsmath}
\usepackage{graphicx}
\usepackage{amsthm}
\usepackage{indentfirst}
\input xy
\xyoption{all}

\newtheorem{defn}{Definition}[section]
\newtheorem{theo}[defn]{Theorem}
\newtheorem{prop}[defn]{Proposition}
\newtheorem{rem}[defn]{Remark}
\newtheorem{tbl}[defn]{Table}
\newenvironment{dem}{\textbf{Proof.}\;}{\hspace{\stretch{1}}\rule{1ex}{1ex}}

\newcommand{\spa}{\mathbb{P}^2(\mathbb{C})}
\newcommand{\lin}{\mathbb{P}^1(\mathbb{C})}
\newcommand{\pto}{\mathbb{P}^0(\mathbb{C})}
\newcommand{\ssm}{X^{SS}(m)}
\newcommand{\sm}{X^S(m)}
\newcommand{\cq}{X^{SS}(m)//G}
\newcommand{\gq}{X^S(m)/G}
\newcommand{\gc}{\gamma^C}
\newcommand{\gl}{\gamma^L}
\newcommand{\Ul}{U^L}
\newcommand{\Uc}{U^C}
\newcommand{\mw}{\widehat{m}}
\begin{document}

\title{GIT quotients of products of projective planes}

\author{ Francesca Incensi }

\date{January 2008.}

\maketitle

\begin{abstract}
We study the quotients for the diagonal
action of $SL_3(\mathbb{C})$ on the product of $n$-fold of $\spa$:
we are interested in describing how the quotient changes when we
vary the polarization (i.e. the choice of an ample linearized line
bundle). We illustrate the different techniques for the construction
of a quotient, in particular the numerical criterion for
semi-stability and the ``elementary transformations'' which are
resolutions of precisely described singularities (case $n=6$).
\end{abstract}


\section*{Introduction}

\indent Consider a projective algebraic variety $X$ acted on by a
reductive algebraic group $G$. Geometric
Invariant Theory (GIT) gives a construction of a $G$-invariant open
subset $U$ of $X$ for which the quotient $U//G$ exists and $U$ is
maximal with this property (roughly speaking, $U$ is obtained by $X$
throwing away ``bad'' orbits). However the open $G$-invariant subset $U$
depends on the choice of a
$G$-linearized ample line bundle. Given an ample $G$-linearized line
bundle $L \in \textrm{Pic}^G(X)$ over $X$, one defines the set of
semi-stable points as 
\[X^{SS}(L):=\{x \in X\,|\, \exists n >0 
\textrm{ and } s \in \Gamma(X,L^{\otimes n})^G\, \textrm{s.t. } s(x)\neq 0\},\] 
and the set of stable points as 
\[X^S(L):=\{x \in X^{SS}(L)\,|\, G\cdot x \textrm{ is closed in } X^{SS}(L) 
\textrm{ and  the stabilizer } G_x \textrm{ is finite} \}.\] 
Then it is possible to introduce a categorical quotient $X^{SS}(L)//G$ in which
two points are identified if the closure of their orbits intersect.
Moreover as shown in \cite{MuFoKi}, $X^{SS}(L)//G$ exists as a
projective variety and contains the \emph{orbit space} $X^S(L)/G$ as
a Zariski open subset.
\[
\begin{array}{cccc} X & & & \\ \cup & & & \\ X^{SS}(L) &
\stackrel{\phi}{\longrightarrow} & X^{SS}(L)//G\\ \cup & & \cup&
\\ X^{S}(L) &
\stackrel{\phi_{|X^S(L)}}{\longrightarrow} & X^{S}(L)/G
\end{array}
\]
{\bf{Question.}} If one fixes $X, G$ and the action of $G$ on $X$,
but lets the linearized ample line bundle $L$ vary in
$\textrm{Pic}^G(X)$, how do the open
set $X^{SS}(L)\subset X$ and the quotient $X^{SS}(L)//G$ change?\\
Dolgachev-Hu \cite{DoHu} and Thaddeus \cite{Tha} proved that only a
finite number of GIT quotients can be obtained when $L$ varies and
gave a general description of the maps relating the various
quotients.
\\\\\indent In this paper we study the geometry of the GIT
quotients for $X=\spa\times\ldots\times \spa=\spa^n$. We give
examples for $n=5$ and $n=6$. The contents of the paper are more
precisely as follows.\\ \indent Section 1 treats the general case
$X=\spa^n$: first of all the numerical criterion of semi-stability
is proved (Proposition 1.1). By means of this it is possible to show
that only a finite number of quotients $\cq$ exists (Subsection
1.2). At the end of the section we introduce the elementary transformations
which relate the different quotients.
\\ \indent Section 2 is concerned with
the case $n=5$. Theorem \ref{theo:p25} contains the main result of
Section 2: we show that there are precisely six different quotients.\\
\indent Section 3 discusses the case $n=6$: the main results of this
Section are concerned with the number of different geometric
quotients that may be obtained (it is 38: Table \ref{table1}) and
with the singularities that may appear in the quotients. In
particular there are only two different types of singularities: in
Subsection 3.2 they are described, using the \'{E}tale
Slice theorem. Theorem \ref{theo:p26} collects these results. At the
end of the Section two examples shows how these
singularities are resolved by ``crossing the wall''.


\subsection*{Acknowledgment}
The results of this paper were obtained during my Ph.D. studies at
University of Bologna and are also contained in my thesis \cite{Inc}
with the same title. I would like to express deep gratitude to my
supervisor prof. Luca Migliorini, whose guidance and support were
crucial for the successful completion of this project.


\section{The general case $X=\spa^n$}
Let $G$ be the group $SL_3(\mathbb{C})$ acting on the variety
$X=\spa^n$ and let $\sigma$ be the diagonal action
\[\begin{array}{cccrcc}\sigma :&
G&
\times& \spa^n & \rightarrow & \spa^n\\
& \,g&,&(x_1,\ldots,x_n)& \mapsto & (gx_1,\ldots,gx_n)\end{array}\]
A line bundle $L$ over $X$ is determined by $L=L(m):=L(m_1,\ldots,
m_n)=\bigotimes_{i=1}^n \pi_i^*(\mathcal{O}_{\spa}(m_i)),\; m_i \in
\mathbb{Z}\; \forall i\;$, where $\pi_i:X\rightarrow\spa$ is the
$i$-th projection. In particular $L$ is ample iff $m_i>0\,,\;
\forall i$.\\ Moreover since each $\pi_i$ is an $G$-equivariant
morphism, $L$ admits a canonical $G$-linearization:
\[\textrm{Pic}^G(X)\cong \mathbb{Z}^n\,.\] 
Thus a \emph{polarization} is completely determined by the line bundle
$L$.
\\\indent Recall that a point $x\in X$ is said to be
\emph{semi-stable} with respect to the polarization $m$ iff there
exists a $G$-invariant section of some positive tensor power of $L$,
$\gamma \in \Gamma(X,L^{\otimes k})^G$, such that $\gamma(x)\neq 0$.
A semi-stable point is \emph{stable} if its orbit is closed and has
maximal dimension. The \emph{categorical quotient} of the open set
of semi-stable points exists and is denoted by $\cq$:
\[\cq\,\cong\, \textrm{Proj}\left(\bigoplus_{k=0}^\infty \Gamma(X,L^{\otimes k})^G\right)\,.\]
Moreover the open set $\gq$ of $\cq$ is a \emph{geometric
quotient}.\\ We set $X^{US}(m)=X\setminus \ssm$, the closed set
of unstable points and $X^{SSS}(m)=\ssm\setminus\sm$, the set of
strictly semi-stable points.


\subsection{Numerical Criterion of semi-stability}
Fixed a polarization $L(m)$, we want to describe the set of
semi-stable points $\ssm$: using the Hilbert-Mumford numerical
criterion, we prove the following
\begin{prop} Let $x\in X$. Then we have
\begin{equation}\label{eq:sss} x\in \ssm \Leftrightarrow
\left\{
\begin{array}{l}\sum_{k, x_k = y}m_k \leq \frac{|m|}{3}\\\\
\sum_{j, x_j \in r} m_j \leq 2\frac{|m|}{3}
\end{array}\right.\end{equation} where $|m|:=\sum_{i =1}^n m_i\,$, and  $y, r$
are respectively a point and a line in $\spa$.
\end{prop}
\begin{dem} Fixing projective coordinates on the $i$-th copy of $\spa$,
$[x_{i0}:x_{i1}:x_{i2}]$, a point $x\in X \left(\subset
\mathbb{P}(\Gamma(X,L(m))^*)=\mathbb{P}^N(\mathbb{C})\right)$, is
described by homogeneous coordinates of this kind:
\[\prod_{i=1}^n x_{i0}^{j_i}x_{i1}^{k_i}x_{i2}^{m_i-(j_i+k_i)}\]
where $0\leq j_i,k_i\leq m_i$, $j_i+k_i\leq m_i$.
\\\indent Let
$\lambda_{\alpha_0,\alpha_1,\alpha_2}$ a one-parameter subgroup of
$G$; it is defined by $
\lambda_{\alpha_0,\alpha_1,\alpha_2}(t)=\textrm{diag}(t^{\alpha_0},t^{\alpha_1},t^{\alpha_2})$
where $\alpha_0+\alpha_1+\alpha_2=0$; we can assume $\alpha_0\geq
\alpha_1\geq \alpha_2$.\\ The subgroup
$\lambda_{\alpha_0,\alpha_1,\alpha_2}$ acts on every component of
$\mathbb{C}^{N+1}$, multiplying by
\[t^{\alpha_0\sum_i j_i+ \alpha_1\sum_i
k_i+\alpha_2\sum_i(m_i-(j_i+k_i))}\,.\] By the definition of the
numerical function of Hilbert-Mumford $\mu_L(x,\lambda)$, we are
interested in determining the minimum value of
\[\alpha_0\sum_{i=1}^n j_i+ \alpha_1\sum_{i=1}^n
k_i+\alpha_2\sum_{i=1}^n\left(m_i-(j_i+k_i)\right)\,.\] This should
be obtained when $j_i=k_i=0$,\, $\forall i$\,; but if there are some
$x_{i2}=0$, then the minimum value becomes:
\begin{equation}\label{eq:mininum} \alpha_2\sum_{i, x_{i2}\neq 0}m_i+
\:\alpha_{1}\sum_{j, x_{j2}=0,
x_{j1}\neq 0} m_j+\:\alpha_0\sum_{k, x_{k2}=x_{k1}=0}m_k\,.
\end{equation}
Thus $x\in X$ is semi-stable if and only if expression
(\ref{eq:mininum})
is less or equal than zero.\\
Let
\[\alpha_0=\beta_0+\beta_1\,,\quad \alpha_1=-\beta_0\,,\quad
\alpha_2=-\beta_1\,;\] it follows that
$\beta_1\geq-2\beta_0$, $\beta_1\geq \beta_0$ e $\beta_1\geq 0$.\\
The expression (\ref{eq:mininum}) can be rewritten and the minimum
value is
\begin{equation}\label{beta}\beta_0\left(\sum_{k, x_{k2}=x_{k1}=0}m_k-\sum_{j,
x_{j2}=0, x_{j1}\neq 0} m_j \right)+\beta_1\left(\sum_{k,
x_{k2}=x_{k1}=0}m_k-\sum_{i, x_{i2}\neq 0}m_i\right)\leq
0\end{equation}
\begin{figure}[h]
\begin{center}
\includegraphics[bb=-70 320 615 780,
scale=0.6]{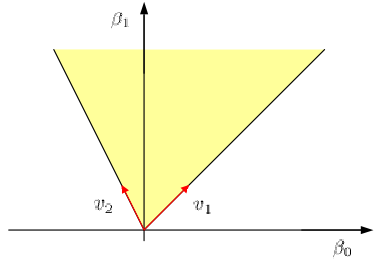} \vspace{-7.2cm}\caption{Plane $\beta_0,
\beta_1$}\label{fig:plane}
\end{center}
\end{figure}
The figure \ref{fig:plane} shows that every couple ($\beta_0,
\beta_1$) that satisfies (\ref{beta}) is a positive linear
combination of $v_1=(1,1)$ e $v_2=(-1,2)$. Thus the relation
(\ref{beta}) must be verified in the two cases $\beta_0=\beta_1=1$ e
$\beta_0=-1, \beta_1=2$.
After few calculations we obtain
\[\left\{\begin{array}{l}\sum_{h, x_h=y}m_h\,\leq
|m|/3\,,\quad y \in \spa\,;\\
\sum_{l, x_l\in r}m_l\,\leq 2|m|/3\,,\quad r \subset
\spa\,.\end{array}\right.\] \end{dem}
\begin{rem} $x \in \sm$ iff the numerical criterion (\ref{eq:sss})
is verified with strict inequalities.\end{rem}
\indent The numerical criterion can be restated as follows: if $K,J$
are subset of $[n]:=\{1,\ldots,n\}$, then we can associate them with
the numbers:
\[\gc_K(m)\,= |m|-3 \sum_{k \in K} m_k \,,\qquad
\gl_J(m)\, = 2|m|-3 \sum_{j \in J} m_j \,. \] In particular we have:
$ \gc_J(m)\,=\, -\gl_{J'}(m)$\; where $ J'=[n]\setminus J$.\\
\indent Now for every subset $K\subseteq [n]$, we consider the set
of configurations $(x_1,\ldots,x_n)$ where the points indexed by $K$
are coincident, while the others are all distinct:
\[\Uc_K=\{x \in X\,|\, x_{k_1}=\ldots=x_{k_{|K|}} \neq x_i ,\, x_j
\neq x_l \: \forall i,j,l \notin K \}\,;\] if $U_K^C \subset \ssm$,
then $\gc_K(m)\geq 0$\,.\\
In the same way if $r$ is a fixed line of $\spa$, let
\[\Ul_J=\{x \in X\,|\, x_{j_1},.., x_{j_{|J|}} \in r, x_i \notin r,
x_i, x_k, x_l\: \textrm{not collinear}\,,\; \forall i,k,l \notin J
\},\] the set of configurations $(x_1,\ldots,x_n)$ where the points
indexed by $J$ are collinear, while the others are not; if $\Ul_J
\subset \ssm$\,, then $\gl_J(m)\geq 0$\,.


\subsection{Quotients}
\begin{prop}\label{prop:GEN}
Let
\[U^{\textrm{GEN}}:=\{x \in X |\; x_1,\ldots, x_n \,
\textrm{ in general position} \}\,\subset X\,,\]
(i.e. every four points among $\{x_1,\ldots, x_n\}$ are a
projective system of $\spa$).\\
Then:
\begin{enumerate}
\item $\ssm\neq \emptyset\, \Leftrightarrow\, U^{\textrm{GEN}} \subset
\ssm$;
\item $\sm \neq \emptyset\, \Leftrightarrow\, U^{\textrm{GEN}} \subset \sm\,
\Leftrightarrow\,\dim(\cq) = 2(n-4)$.
\end{enumerate}
\end{prop}
We know that the quotient $\cq$ depends on the choice of the
polarization $L(m)$: moreover Dolgachev-Hu \cite{DoHu} and Thaddeus
\cite{Tha} have proved that when $L(m)$ varies, then there exists
only a \textsl{finite} number of different quotients.\\\indent Now
we give a proof of the same result in our case.\\ If $\ssm\neq
\emptyset$, then by the previous Proposition we have
$U^{\textrm{GEN}}\subset \ssm$. Moreover sets $\Uc_K$ and $\Ul_J$
are in a finite number since
they consist in particular combinations of $x_1,\ldots, x_n$.\\
Fixed a polarization $m$, $\ssm$ can be described as
\[\ssm=U^{\textrm{GEN}}\,\cup\, {\cal{U^{SS}}}(m)\,,\] where
${\cal{U^{SS}}}(m):=\left\{\,\Uc_K,
\Ul_J\,|\, \Uc_K, \Ul_J \subset \ssm\,\right\}$. In particular we
can construct only a finite number of different sets
${\cal{U^{SS}}}(m)$ and as a consequence there exists a finite
number of different open sets $\ssm$; in conclusion only a finite
number of quotients $\cq$ exists.


\subsection{Elementary transformations}
Let $m$ be a polarization such that $3$ divides $|m|$ and $\sm\neq
\emptyset, \sm \subsetneq \ssm$; let us consider ``variations'' of
$m$ as follows:
\[\mw = m \pm
(0,\ldots,0,\underbrace{1}_{i},0,\ldots,0)\,.\] We can have two
different kind of variations, depending on the value $|\mw|$\,:
\begin{enumerate}
\item $\mw\,
\stackrel{+1_i}{\longrightarrow}\, m$ \, (i.e. $|\mw| \equiv 2\quad
\textrm{mod }3$)\,;
\item $\mw\,\stackrel{-1_i}{\longrightarrow}\,
m$ \, (i.e. $|\mw| \equiv 1\quad \textrm{mod }3$)\,.\\
\end{enumerate}
In both cases we have $X^S(\mw)=X^{SS}(\mw)$; studying the relations
between values $\gc_J(\mw), \gl_K(\mw)$ and values $\gc_J(m),
\gl_K(m)$, we observe that
\begin{itemize}
\item[\textbf{1.}] \quad $\mw
\stackrel{+1_i}{\longrightarrow} m$\;
\[X^S(\mw) \subset \ssm\,,\quad
X^S(\mw)= \ssm \setminus \bigcup_{i \notin J, \gc_J(m)=0\, \vee\,
\gl_J(m)=0} U_J^*\;; \]
\[\sm\subset X^S(\mw)\,,\quad
\sm= X^S(\mw)\setminus \bigcup_{i \in H, \gc_H(\mw)=2\, \vee\,
\gl_H(\mw)=1} U_H^*\,.\]
\item[\textbf{2.}]\quad $\mw
\stackrel{-1_i}{\longrightarrow} m$
\[X^S(\mw) \subset \ssm\,,\quad
X^S(\mw)= \ssm \setminus \bigcup_{i \in J, \gc_J(m)=0\, \vee\,
\gl_J(m)=0} U_J^*\;; \]
\[\sm\subset X^S(\mw)\,,\quad
\sm=X^S(\mw) \setminus \bigcup_{i \notin H, \gc_H(\mw)=1\, \vee\,
\gl_H(\mw)=2} U_H^*\,. \]
\end{itemize}
\noindent At the end, we can illustrate the inclusions of the open
sets of stable and semi-stable points, with the following diagrams:
\[ \xymatrix{
X^{SS}(\widehat{m})=\hspace{-1cm} &X^S(\widehat{m})\ar@{^{(}->}[r]^\alpha & X^{SS}(m) \\
&X^S(m)\ar[u]^\beta \ar[ur] &}\]
\setlength{\unitlength}{1cm}
\[\begin{picture}(5,5)
\linethickness{0.02mm} \put(2.5,3.5){\oval(6.5,3.5)}
\put(2.5,3){\oval(4,2.25)}
 \put(2.5,2.5){\oval(2,1)}
 \put(2,2.3){$\sm$}
\put(2,4.5){$\ssm$}
 \put(2,3.3){$X^S(\mw)$}
\end{picture}
\vspace{-1.5cm}\] 
\indent The inclusions $\sm\subset X^{S}(\mw)\subset
\ssm$ induce a morphism
\begin{equation} \theta: X^S(\mw)/G \longrightarrow
\cq\,, \end{equation} 
which is an isomorphism over $\gq$, while over
$\left(\cq\right)\setminus \left(\gq\right)$ is a contraction of
subvarieties.
\\\indent In fact, let us consider a point $\xi \in
\left(\cq\right)\setminus \left(\gq\right)$: this is the image in
$\cq$ of different open, strictly semi-stable orbits, that all have
in their closure a closed, minimal orbit $Gx$, for a certain
configuration $x=(x_1,\ldots,x_n)\in X^{SSS}(m)$. In particular this
configuration $x$ has $|J|$ coincident points, and the others
$n-|J|$ collinear; by the numerical criterion, we get $\gc_J(m)=0$
and $\gl_{J'}(m)=0$, where $J$ indicates the coincident points,
while $J'=[n]\setminus J$ indicates the collinear ones.\\
For the sake of simplicity, we can assume $x$ as
\[\left(\begin{array}{ccccccccc}
1& \ldots& 1& 0& 0&0& 0&\ldots&0\\
0&\ldots&0&1&0&1&1&\ldots&
1\\0&\ldots&0&0&1&\beta_1&\beta_2&\ldots&\beta_{n-|J|-2}\end{array}
\right)\,,\quad\beta_k \in \mathbb{C}^*,\; \forall k\,.\]
\\\indent The open orbits $O$ that contain $Gx$ in their closure, are
characterized by $\gc_J(m)=0$ or $\gl_{J'}(m)=0$; there are two
different cases:
\begin{enumerate}
\item $\gc_J(m)=0$: orbits look as
\[O_1=\left(\begin{array}{ccccccccc}
1& \ldots& 1& 0& 0&\alpha_1& \alpha_2&\ldots&\alpha_{n-|J|-2}\\
0&\ldots&0&1&0&1&1&\ldots&
1\\0&\ldots&0&0&1&\rho\beta_1&\rho\beta_2&\ldots&\rho\beta_{n-|J|-2}\end{array}
\right), \rho \in \mathbb{C}^*, \alpha_k \in \mathbb{C}.\]
\item $\gl_{J'}(m)=0$: orbits look as
\[\label{orb
all}O_2=\left(\begin{array}{cccccccccc}
1& 1&\ldots& 1& 0& 0&0& 0&\ldots&0\\
0& \delta_1&\ldots&\delta_{|J|-1}&1&0&1&1&\ldots&
1\\0&\epsilon_1&\ldots&\epsilon_{|J|-1}&0&1&\beta_1&\beta_2&\ldots&
\beta_{n-|J|-2}\end{array} \right), \delta_k, \epsilon_k \in
\mathbb{C}.\]
\end{enumerate}
\indent Now, calculating $\theta^{-1}(\xi)$, it follows:
\[\theta^{-1}(\xi)=\theta^{-1}\left(\phi(\overline{\Uc_J} \cup
\overline{\Ul_{J'}})\right)\,;\] by the numerical criterion, only
one between $\Uc_J$ and $\Ul_{J'}$ is included in $X^S(\mw)$.
\\\indent Dealing with an elementary transformation of the first
type ($\mw \stackrel{+1_i}{\longrightarrow} m$), then
\begin{itemize}\item[-]if $i \in J\, \Rightarrow
\theta^{-1}(\xi)=\theta^{-1}\left(\phi(\overline{\Uc_J} \cup
\overline{\Ul_{J'}})\right)= \widehat{\phi}\left(\overline{\Uc_J}
\cap X^S(\mw)\right)\,.$
\\When $n\geq 5$, this has dimension:
\begin{equation}\label{dim fiber i in J}
d=n-|J|-3\,.
\end{equation}
In fact, let us consider the minimal closed orbit $Gx$: all the
orbits that contain $Gx$ in their closure and are stable in
$X^S(\mw)$, are characterized by the coincidence of $|J|$ points
($O_1$ orbits ).
\item[-] if $i \in J'\, \Rightarrow
\theta^{-1}(\xi)=\theta^{-1}\left(\phi(\overline{\Uc_J} \cup
\overline{\Ul_{J'}})\right)=
\widehat{\phi}\left(\overline{\Ul_{J'} }\cap X^S(\mw)\right)\,.$\\
Now the dimension $d$ of $\theta^{-1}(\xi)$ is
\begin{equation}\label{dim fiber i in J'}
 d=2\left(n-|J'|-1\right)-1\,.
\end{equation}
\end{itemize}
\indent Dealing with an elementary transformation of the second type
($\mw \stackrel{-1_i}{\longrightarrow} m$), then
\begin{equation}\label{dim fiber trasf2}
i \in J \Rightarrow  d=2\left(n-|J'|-1\right)-1\,; \qquad
i \in J' \Rightarrow d=n-|J|-3\,.
\end{equation}

\section{$X=\spa^5$}
\subsection{Number of quotients}
\indent Let us study the case $n=5$: $X=\spa^5$. First of all let us
determine how many different quotients we may get when the
polarization varies.\\\indent Let us examine the number of
\emph{Geometric} quotients; let $L(m)$ be a polarization such that
$\sm \neq \emptyset$: by the Proposition \ref{prop:GEN} it follows
that $U^{\textrm{GEN}}\subset \sm$. In particular
\[\sm=U^{\textrm{GEN}}\,\cup\, {\cal{U^S}}(m)\,,\] where
${\cal{U^S}}(m):=\left\{\,\Uc_K, \Ul_J\,|\, \Uc_K, \Ul_J \subset
\sm\,\right\}$. Obviously there is only a finite number of sets
${\cal{U^{S}}}(m)$: we want to describe their structure.\\
\indent Let $m=(m_1,\ldots,m_5)$ be a polarization such that
$\sm=\ssm\neq \emptyset$; we can assume $m_i\in\mathbb{Q}$ and
\[0<m_i<\frac{1}{3}\,,\quad m_i\geq m_{i+1}\,,\quad |m|=1\,.\] 
As a consequence only strictly inequalities are allowed
in the numerical criterion: \[x \in \sm
\Leftrightarrow\, \sum_{k, x_k=y, k\in K}m_k < \frac{1}{3}\,,
\sum_{j, x_j \in r, j\in J} m_j < \frac{2}{3}\,\Leftrightarrow\,
\gc_K(m)>0\,, \gl_J(m)>0\] In particular sets $K$ that indicate
coincident points, can have only two elements (otherwise it would be
possible to find a weight $m_i$ greater than $1/3$), and in the same
way non trivial sets $J$ that indicate collinear points, have only
three elements.\\ Moreover by the numerical criterion, only some
sets $\Uc_K, \Ul_J$ may be included in $\sm$:
\begin{equation}\label{sets:5}\begin{array}{lllll}U^C_{15}\,,&
\Uc_{25}\,,&\Uc_{34}\,,& \Uc_{35}\,,&\Uc_{45}\,,\\
\Ul_{234}\,,& \Ul_{134}\,,&
\Ul_{125}\,,&\Ul_{124}\,,&\Ul_{123}\,.\end{array} \end{equation}
They can be examined in couple, because $\gc_K(m)=-\gl_{K'}(m)\,,
K'=[5]\setminus K$ and then only one between $\Uc_K$ and $\Ul_{K'}$ may be included in $\sm$.\\
The number of geometric quotient is \emph{six}.
\\In fact
\begin{enumerate}
\item[0.] in ${\cal{U^S}}(m)$ there may be only sets as $\Ul_{J}$: an
example is the polarization $m=\left(1/5, 1/5, 1/5, 1/5,
1/5\right)$\,;
\item if in ${\cal{U^S}}(m)$ there is one set as $\Uc_K$, it is $\Uc_{45}$: in fact
if it were $\Uc_{34}$, then it follows \[m_3+m_4< 1/3\; \textrm{ and
}\; m_4+m_5> 1/3\,\Rightarrow\, m_3<m_5\,\Rightarrow\,
\textrm{Impossible}. \] Example: $m=(1/4,1/4,1/4,1/8,1/8)$\,;
\item if in ${\cal{U^S}}(m)$ there are two sets as $\Uc_K$, they are $\Uc_{45}$ and $\Uc_{35}$:
 the argument is similar to the previous one.\\ Example:
$m=(3/11,3/11,2/11,2/11,1/11)$\,;
\item if in ${\cal{U^S}}(m)$ there are three sets as $\Uc_K$, we can have two cases:
\begin{enumerate}
\item $\Uc_{45}, \Uc_{35}$ and $\Uc_{25}$\,, example
$m=(3/10,1/5,1/5,1/5,1/10)$\,;
\item $\Uc_{45}, \Uc_{35}$ and $\Uc_{34}$\,,  example
$m=(3/10,3/10,1/5,1/10,1/10)$\,.
\end{enumerate}
\item if in ${\cal{U^S}}(m)$ there are four sets as $\Uc_K$,
they are $\Uc_{45}, \Uc_{35}, \Uc_{25}$ and $\Uc_{15}$. Example:
$m=(1/4,1/4,1/4,2/9,1/36)$\,;
\item the case of all $\Uc_K$ sets in ${\cal{U^S}}(m)$ is impossible,
because $\Uc_{45},$ $\Uc_{35},$ $\Uc_{34},$ $\Uc_{25}$ are
incompatible.
\end{enumerate}
We have found six cases:
\[\begin{array}{ll}0.&{\cal{U^S}}(m)=\{\Ul_{234},\Ul_{134},\Ul_{124},\Ul_{123},\Ul_{125}\}\\
1.&{\cal{U^S}}(m)=\{\Ul_{234},\Ul_{134},\Ul_{124},\Ul_{125},\Uc_{45}\}\,,\\
2.&{\cal{U^S}}(m)=\{\Ul_{234},\Ul_{134},\Ul_{125},\Uc_{35},\Uc_{45}\}\,,\\
3a.&{\cal{U^S}}(m)=\{\Ul_{234},\Ul_{125},\Uc_{25},\Uc_{35},\Uc_{45}\}\,,\\
3b.&{\cal{U^S}}(m)=\{\Ul_{234},\Ul_{134},\Uc_{34},\Uc_{35},\Uc_{45}\}\,,\\
4.&{\cal{U^S}}(m)=\{\Ul_{125},\Uc_{15},\Uc_{25},\Uc_{35},\Uc_{45}\}\,.\end{array}\]
Then there are only six different open sets of stable points and
thus six geometric quotients.\\
\indent Now let us examine the number of \emph{Categorical}
quotients. First of all let us observe that sets $\Uc_K, \Ul_{K'}$
that may be included in $\ssm$ are the same of (\ref{sets:5}). What
is different from the previous case is that now two sets $\Uc_K$ and
$\Ul_{K'}$ may be \emph{both} included in $\ssm$ (if
$\gc_K(m)=\gl_{K'}(m)=0$); this means that in $\ssm$ there are two
distinct strictly semi-stable orbits:
\begin{itemize}
\item[-] an orbit $O_1$ with $x_{k_1}=x_{k_2},\, K=\{k_1,k_2\}$;
\item[-] orbits $O_2$ with $x_{i_1}, x_{i_2}, x_{i_3}$ collinear, $i_1, i_2, i_3 \in K'$.
\end{itemize}
Orbit $O_1$ and all orbits $O_2$ contain in their closure a closed,
minimal, strictly semi-stable orbit $O_{12}$, that is characterized
by $x_{k_1}=x_{k_2}$ \emph{and} $x_{i_1}, x_{i_2}, x_{i_3}$
collinear: \setlength{\unitlength}{1cm}
\begin{picture}(4,4)
 \put(4,1){\line(2,1){4.5}}
 \put(5,1.5){\circle*{0.15}}
 \put(6,2){\circle*{0.15}}
 \put(7,2.5){\circle*{0.15}}
 \put(5,3){\circle{0.15}}
 \put(4.3,3.2){$x_{k_1}=x_{k_2}$}
 \put(5.1,1.1){$x_{i_1}$}
 \put(6.1,1.6){$x_{i_2}$}
 \put(7.1,2.1){$x_{i_3}$}
\end{picture}\vspace{-.5cm}\newline
In the categorical quotient $\cq$, orbits $O_1$ and $O_2$ determine
the \emph{same} point; in fact $O_{12}\subset(\overline{O_1} \cap
\overline{O_2})$. \\\indent Let us examine the stable case more
accurately: we know that only one between $O_1$ and $O_2$ is
included in $\sm$; when $O_1$ is included, it determines a point of
the geometric quotient. In fact if for example $\Uc_{45}\subset
\sm$, then $\phi(\Uc_{45})$ may regarded as
$\spa^4(m_1,m_2,m_3,m_4+m_5)/SL_3(\mathbb{C})$ and the we have a
point. When orbits $O_2$ are included in $\sm$, they determine a
$\lin$ in $\gq$. In fact if for example $\Ul_{123}\subset \sm$, then
we can assume \[O_2=\left(\begin{array}{ccccc} 1& 0 & 1 & 0 & \alpha\\ 0& 1 & 1 & 0 & \beta\\
0 & 0 & 0 & 1 & 1 \end{array}\right), (\alpha,\beta)\in
\mathbb{C}^2\setminus\{(0,0)\}\,.\] Applying to $O_2$ a projectivity
$G_{\lambda}$ of $\spa$ that fixes the line that contains $x_1, x_2,
x_3$ ($G_{\lambda}\cong\textrm{diag}(\lambda,\lambda,\lambda^{-2}),$
with $\lambda \in \mathbb{C}^*$), , it follows:
\[G_\lambda\cdot x=\left(\begin{array}{ccccc} 1& 0 & 1 & 0 & \lambda^3 \alpha\\
0& 1 & 1 & 0 & \lambda^3\beta\\
0 & 0 & 0 & 1& 1
\end{array}\right)\,.\]
If $\alpha\neq 0$, then we can assume $\lambda^3=\alpha^{-1}$; thus
we obtain $x_5=[1:\alpha^{-1}\beta:1]$; in the same way if
$\beta\neq 0$, then $x_5=[\alpha\beta^{-1}:1:1]$. \\ Then it is
clear that $\phi(O_2)\cong\lin$.\\\indent In the semi-stable case
when $\Uc_K, \Ul_{K'} \subset \ssm$, we know that $\overline{\Uc_K}
\cap \overline{\Ul_{K'}} \neq \emptyset$ and they determine a
non-singular point of $\cq$, just as in the stable case when
$\Uc_K\subset\sm$. In this way it follows that every categorical
quotient $\cq$, where \[\begin{array}{ccc}\ssm=U^{\textrm{GEN}}\cup
\{&\underbrace{\Uc_J,\Ul_I,\ldots,}&
\underbrace{\Uc_K,\Ul_{K'},\ldots,\Uc_{H},\Ul_{H'}}\},\\ &
 stable sets & semi-stable sets\end{array}\] is isomorphic to a geometric one $X^S(m')/G$, whose
open set of stable points is
\[X^S(m')=U^{\textrm{GEN}}\cup\{\Uc_J,\Ul_I,\ldots,
\Uc_K,\Uc_{\cdots},\Uc_{H},\}.\]
\\\indent In conclusion:
\begin{theo}\label{theo:p25} Let $X=\spa^5$: then there are \emph{six}
non trivial quotients.\\Moreover a quotient $\cq$ is isomorphic to
one of the following:
\[\begin{array}{ll}\spa &\\
\spa \textrm{with a point blown up} &(\spa_1)\\
\lin\times\lin& (\lin^2);\\
\spa\textrm{with two points blown up}& (\spa_2)\\
\spa \textrm{with three points blown up}& (\spa_3)\\
\spa \textrm{with four points blown up}& (\spa_4) \end{array}\]
\end{theo}

\subsection{Quotients $\spa^5//G$}
The following diagram shows the relations between some polarizations that realize the quotients; for example if $m=(22211)$, then $X^S(m)=\spa_3$ and there is a morphism $\theta:X^S(22211)/G=\spa_3\rightarrow X^{SS}(44322)//G=\lin\times\lin$.
\[ \xymatrix{
&&(11111)_{\spa_4}\ar@<-2.0ex>[d] &\\
&&(21111)_{\lin} &\\
&&(22111)_{(\lin)^2}\ar@<2.0ex>[u]\ar[dl]& \\
&(44322)_{(\lin)^2}&&\\
&&(22211)_{\spa_3}\ar[ul]\ar[dl]\ar[dr]&\\
(33111)_{\pto}& (32211)_{\lin} &&(22221)_{\spa}\\
&(33211)_{(\lin)^2}\ar[ul]\ar@<2.0ex>[u]\ar[dl]\ar[d]&(32221)_{\spa_1}\ar[ul]\ar[ur]\ar@<-2.0ex>[d]&(22222)_{\spa_4}\ar@<2.0ex>[u]\ar@<-2.0ex>[d]\\
(66522)_{(\lin)^2}&(66432)_{(\lin)^2}&(65442)_{\spa_1}&(54444)_{\spa_4}\\
&(33311)_{\spa_3}\ar[ul]\ar@<-2.0ex>[d]&(33221)_{\spa_2}\ar[ul]\ar@<2.0ex>[u]\ar[dr]\ar[dl]&(32222)_{\spa_4}\ar@<2.0ex>[u]\ar@<-2.0ex>[d]\\
&(33321)_{\spa}&&(33222)_{(\lin)^2}\\
&&(33322)_{\spa_3}\ar[ul]\ar[ur]&\\
} \]

\section{$X=\spa^6$}
\subsection{Number of quotients}
Now we study the case $n=6$: $X=\spa^6$; as in the previous case we
first determine how many different \emph{Geometric} quotient we can
get when the polarization varies.\\\indent For a polarization
$m=(m_1,\ldots,m_6)$ such that $\sm\neq \emptyset$, then
\[\sm=U^{\textrm{GEN}}\,\cup\, {\cal{U^S}}(m)\,.\]
We want to describe the structure of the sets ${\cal{U^{S}}}(m)$;
assume that $0<m_i<\frac{1}{3}$\,, $m_i\geq m_{i+1}$\,, $|m|=1\,.$\\
We are interested in those sets $\Uc_K$ that are included in $\sm$:
some are \textsl{always} included in $\sm$:
\[\Uc_{36}\,,\quad \Uc_{46}\,,\quad \Uc_{56}\,,\] and others \textsl{may} be
included in $\sm$:
\[\begin{array}{lllllllll} \Uc_{15}\,,&
\Uc_{16}\,,&\Uc_{23}\,,& \Uc_{24}\,,& \Uc_{25}\,,& \Uc_{26}\,,&
\Uc_{34}\,,& \Uc_{35}\,, & \Uc_{45}\,,\\ \Uc_{156}\,,& \Uc_{256}\,,&
\Uc_{345}\,,&\Uc_{346}\,,&\Uc_{356}\,& \Uc_{456}\,.\end{array}
\]
The number of different sets ${\cal{U^{S}}}(m)$ is 38.\\\indent
First of all the minimum number of sets $\Uc_K$ with $|K|=2$,
included in $\sm$ is five: in fact for example consider only the
sets $\Uc_{36}, \Uc_{46}, \Uc_{56}$ that are always included in
$\sm$, then obviously \[m_1+m_6>\frac{1}{3}\,,\,
m_2+m_5>\frac{1}{3}\,,\, m_3+m_4>\frac{1}{3}\; \Rightarrow\,
\sum_{i=1}^6 m_i >1\,:\; \textrm{impossible}.\] In a similar way it
is impossible to have only four sets $\Uc_K\, (|K|=2)$ in $\sm$.\\
\indent Then for five sets $\Uc_K$, we have $\Uc_{16}, \Uc_{26},
\Uc_{36}, \Uc_{46}, \Uc_{56}$: in fact with another $5$-tuple (for
example $\Uc_{45}, \Uc_{26}, \Uc_{36}, \Uc_{46}, \Uc_{56}$), it gets
$|m|>1\,,$  that is impossible. \\ Moreover with these combinations,
it is impossible to obtain a set as $\Uc_K$ with $|K|=3$.\\ \indent
Going on with the calculations, we are able to construct the
following table, that shows all the possible cases (in the ``admissible'' cells we exhibit an example of a polarization that realize the geometric quotient). In particular it
is not possible to have more than ten sets $\Uc_K\, (|K|=2)$ in
$\sm$: we would obtain $|m|<1$.\\
\begin{tbl}\label{table1}
\[\small{\begin{tabular}{|c|c|c|c|c|c|}
  \hline
  & & & & & \\
  $\Uc_K, $ & No $\Uc_K$, & $1$ set $\Uc_K, $ & $2$ sets $\Uc_K, $ & $3$ sets $\Uc_K, $  & $4$ sets $\Uc_K, $ \\
  $|K|=2$&  $|K|=3$ & $|K|=3$& $|K|=3$& $|K|=3$&$|K|=3$ \\
  & & & & & \\
  \hline\hline
  & & & & & \\
  $\Uc_{16},\Uc_{26},\Uc_{36},$ & \checkmark & No$^{(*)}$ & No & No & No \\
  $\Uc_{46},\Uc_{56}$ & $
  \frac{1}{11}(222221)$ & &  & & \\
  & & & & & \\
  \hline
  & & & & & \\
  $\Uc_{16},\Uc_{26},\Uc_{36},$ & \checkmark & $\Uc_{456}$ & No$^{(*)}$ & No & No \\
  $\Uc_{45},\Uc_{46},\Uc_{56}$ & $\frac{1}{14}(333221)$
  & $\frac{1}{17}(444221)$ & & & \\
  & & & & & \\
  \hline
  & & & & & \\
  $\Uc_{34},\Uc_{35},\Uc_{36},$ & \checkmark & $\Uc_{456}$ & $\Uc_{456},\Uc_{356}$ &
  $\Uc_{456},\Uc_{356},$& $\Uc_{456},\Uc_{356},$ \\
  $\Uc_{45},\Uc_{46},\Uc_{56}$ & & & & $\Uc_{346}$ & $\Uc_{346},\Uc_{345}$ \\
  & $\frac{1}{8}(221111)$ & $\frac{1}{11}(332111)$ & $\frac{1}{14}(442211)$ &
  $\frac{1}{17}(552221)$ & $\frac{1}{10}(331111)$ \\
  & & & & & \\
  \hline
  & & & & & \\
  $\Uc_{25},\Uc_{26},\Uc_{35},$ & \checkmark & $\Uc_{456}$ & $\Uc_{456},\Uc_{356}$ & $\Uc_{456},\Uc_{356},$ & No$^{(*)}$ \\
  $\Uc_{36},\Uc_{45},\Uc_{46},$ & & & & $\Uc_{256}$ & \\
  $\Uc_{56}$ & $\frac{1}{11}(322211)$ & $\frac{1}{14}(433211)$ & $\frac{1}{17}(543311)$ & $\frac{1}{19}(644311)$ & \\
  & & & & & \\
  \hline
  & & & & & \\
  $\Uc_{26},\Uc_{34},\Uc_{35},$ & \checkmark & $\Uc_{456}$ & $\Uc_{456},\Uc_{356}$ & $\Uc_{456},\Uc_{356},$ & No$^{(**)}$ \\
  $\Uc_{36},\Uc_{45},\Uc_{46},$ & & & & $\Uc_{346}$ & \\
  $\Uc_{56}$ & $\frac{1}{14}(432221)$ & $\frac{1}{17}(543221)$ & $\frac{1}{26}(875321)$ & $\frac{1}{16}(542221)$ & \\
  & & & & & \\
  \hline
  & & & & & \\
  $\Uc_{16},\Uc_{26},\Uc_{35},$ & \checkmark & $\Uc_{456}$ & $\Uc_{456},\Uc_{356}$ & No$^{(*)}$ & No \\
  $\Uc_{36},\Uc_{45},\Uc_{46},$ & & & & & \\
  $\Uc_{56}$ & $\frac{1}{17}(443321)$ & $\frac{1}{20}(554321)$ & $\frac{1}{26}(775421)$ & & \\
  & & & & & \\
  \hline
  & & & & & \\
  $\Uc_{16},\Uc_{26},\Uc_{34},$ & \checkmark & $\Uc_{456}$ & $\Uc_{456},\Uc_{356}$ & $\Uc_{456},\Uc_{356},$ & No$^{(**)}$ \\
  $\Uc_{35},\Uc_{36},\Uc_{45},$ & & & & $\Uc_{346}$ & \\
  $\Uc_{46},\Uc_{56}$ & $\frac{1}{13}(332221)$ & $\frac{1}{16}(443221)$ & $\frac{1}{19}(553321)$ & $\frac{1}{25}(774331)$ & \\
  & & & & & \\
  \hline
    & & & & & \\
  $\Uc_{16},\Uc_{25},\Uc_{26},$ & \checkmark & $\Uc_{456}$ & $\Uc_{456},\Uc_{356}$ & $\Uc_{456},\Uc_{356},$ & No$^{(*)}$ \\
  $\Uc_{35},\Uc_{36},\Uc_{45},$ & & & & $\Uc_{256}$ & \\
  $\Uc_{46},\Uc_{56}$ & $\frac{1}{16}(433321)$ & $\frac{1}{26}(766421)$ & $\frac{1}{26}(765521)$ & $\frac{1}{25}(755521)$ & \\
  & & & & & \\
  \hline
  & & & & & \\
  $\Uc_{25},\Uc_{26},\Uc_{34},$ & \checkmark & $\Uc_{456}$ & $\Uc_{456},\Uc_{356}$ & No$^{(\dag)}$ & No \\
  $\Uc_{35},\Uc_{36},\Uc_{45},$ & & & & & \\
  $\Uc_{46},\Uc_{56}$ & $\frac{1}{31}(965542)$ & $\frac{1}{26}(865322)$ & $\frac{1}{13}(432211)$ & & \\
  & & & & & \\
  \hline

  \end{tabular}}\]

  \[\small{\begin{tabular}{|c|c|c|c|c|c|}
  \hline
  & & & & & \\
  $\Uc_K, $ & No $\Uc_K$, & 1 set $\Uc_K, $ & 2 sets $\Uc_K, $ & 3 sets $\Uc_K, $ & 4 sets $\Uc_K, $ \\
  $|K|=2$ & $|K|=3$ & $|K|=3$ & $|K|=3$ & $|K|=3$ & $|K|=3$ \\
  & & & & & \\
  \hline
  \hline
  & & & & & \\
  $\Uc_{15},\Uc_{16},\Uc_{25},$ & \checkmark & $\Uc_{456}$ & $\Uc_{456},\Uc_{356}$ & $\Uc_{456},\Uc_{356},$ & $\Uc_{456},\Uc_{356},$ \\
  $\Uc_{26},\Uc_{35},\Uc_{36},$ & & & & $\Uc_{256}$ & $ \Uc_{256},\Uc_{156}$ \\
  $\Uc_{45},\Uc_{46},\Uc_{56}$ & $\frac{1}{10}(222211)$ & $\frac{1}{13}(333211)$ &
  $\frac{1}{16}(443311)$ & $\frac{1}{25}(766411)$ & $\frac{1}{22}(555511)$ \\
  & & & & & \\
  \hline
  & & & & & \\
  $\Uc_{24},\Uc_{25},\Uc_{26},$ & \checkmark & $\Uc_{456}$ & No$^{(\dag\dag)}$ & No & No \\
  $\Uc_{34},\Uc_{35},\Uc_{36},$ & & & & & \\
  $\Uc_{45},\Uc_{46},\Uc_{56}$ & $\frac{1}{17}(533222)$ & $\frac{1}{10}(322111)$ & & & \\
  & & & & & \\
  \hline
  & & & & & \\
  $\Uc_{23},\Uc_{24},\Uc_{25},$ & \checkmark & No$^{(\dag\dag\dag)}$ & No & No & No \\
  $\Uc_{26},\Uc_{34},\Uc_{35},$ & & & & & \\
  $\Uc_{36},\Uc_{45},\Uc_{46},$ & $\frac{1}{7}(211111)$ & & & & \\
  $\Uc_{56}$ & & & & & \\
  \hline
\end{tabular}}\]
\end{tbl}
\begin{itemize}
\item[$^{(*)}$] This case is not possible, because there is not any available tern;
\item[$^{(**)}$] $\Uc_{345}$ is not included in $\sm$, because otherwise
$m_3+m_4+m_5<\frac{1}{3},$ $ m_2+m_6<\frac{1}{3}$ $\Rightarrow
m_1>\frac{1}{3}$, that is impossible;
\item[$^{(\dag)}$]
$\Uc_{256},\Uc_{345},\Uc_{346}\nsubseteq \sm$;
\item[$^{(\dag\dag)}$]
$\Uc_{246}, \Uc_{256},\Uc_{345},\Uc_{346},\Uc_{356}\nsubseteq \sm$;
\item[$^{(\dag\dag\dag)}$]
$\Uc_{236},\Uc_{246},\Uc_{256},\Uc_{345},\Uc_{346},\Uc_{356},\Uc_{456}\nsubseteq
\sm$.\\
\end{itemize}

\subsection{Singularities}
In this section we study the singularities which appear in the categorical quotients.\\
Suppose that $|m|$ is divisible by 3, and that
there exist strictly semi-stable orbits (included in $X^{SSS}(m)$);
then we can have different cases depending on some ``partitions'' of
the polarization $m \in \mathbb{Z}_{>0}^6$:
\begin{enumerate}
\item  there are two distinct indexes $i,j$ such that $m_i+m_j=|m|/3\,$;
as a consequence, for the other indexes it holds
$m_h+m_k+m_l+m_n=2|m|/3$ (i.e. minimal closed orbits have
$x_i=x_j\,$
and $x_h, x_k, x_l, x_n$ collinear).\\
\setlength{\unitlength}{1cm}
\begin{picture}(4,3.75)
 \put(3,1){\line(3,2){3.5}}
 \put(2.32,3){\circle*{0.15}} \put(1.6,3.2){$x_i=x_j$}
 \put(6,3){\circle{0.15}} \put(6.2,2.8){$x_h$}
 \put(5.25,2.5){\circle{0.15}} \put(5.4,2.3){$x_k$}
 \put(4.4,1.95){\circle{0.15}} \put(4.5,1.7){$x_l$}
 \put(3.7,1.45){\circle{0.15}} \put(3.8,1.2){$x_n$}
 \end{picture}\vspace{-1cm}
\\In $\cq$ these orbits determine a curve $C_{ij}\cong
\lin\,.$
\begin{enumerate}
\item[1.1] particular case: $m_i+m_j=m_h+m_l=m_k+m_n=|m|/3\,$ for distinct indexes
(i.e. there is a ``special'' minimal, closed orbit other than the
orbits previously seen, characterized by
$x_i=x_j\,, x_h=x_l\,, x_k=x_n$\,).\\
\setlength{\unitlength}{1cm}
\begin{picture}(4,3.5)
 \put(2,2.5){\circle*{0.15}} \put(2.1,2.7){$x_i=x_j$}
 \put(4,1.5){\circle*{0.15}} \put(4,1){$x_k=x_n$}
 \put(6,2.5){\circle*{0.15}} \put(6,2.1){$x_h=x_l$}
 \end{picture}\vspace{-1cm}
\end{enumerate}
\item there are three distinct indexes $h,i,j$ such that $m_h+m_i+m_j=|m|/3\,$;
as a consequence for the other indexes it holds $m_k+m_l+m_n=2|m|/3$
(i.e. there is a minimal, closed orbit such that $x_h=x_i=x_j\,$,
and $x_k, x_l, x_n$ collinear).
\\
\setlength{\unitlength}{1cm}
\begin{picture}(6,4)
 \put(4,1){\line(3,2){3}}
 \put(3.32,3){\circle*{0.15}} \put(2.3,3.2){$x_h=x_i=x_j$}
 \put(6.25,2.5){\circle{0.15}} \put(6.4,2.3){$x_k$}
 \put(5.4,1.95){\circle{0.15}} \put(5.5,1.7){$x_l$}
 \put(4.7,1.45){\circle{0.15}} \put(4.8,1.2){$x_n$}
 \end{picture}\vspace{-0.5cm}
\end{enumerate}
\indent Let us study  minimal, closed orbits and what they determine
in $\cq$.

\subsubsection{$x_i=x_j$ and $x_h, x_k, x_l, x_n$ collinear} \indent
Consider a polarization $m=(m_1,\ldots,m_6)$ as previously indicated
and an orbit $Gx$ such that $x_i=x_j\, (m_i+m_j=|m|/3)$, and the
other four points $x_h, x_k, x_l, x_n$ collinear
($m_h+m_k+m_l+m_n=2|m|/3$).
\\$Gx$ is a minimal, closed, strictly semi-stable orbit and its image in
$\cq$ is a point $\xi\in C_{ij}$. For the sake of generality,
suppose that $x_h, x_k, x_l, x_n$ are collinear, but distinct; for
example assume $x$ as:
\[x=\left(
\begin{array}{cccccc} 1 & 1 & 0 & 0 & 0 & 0 \\ 0 & 0 & 1& 0 & 1 &
1 \\ 0 & 0 & 0 & 1 & a& b
\end{array} \right),\quad a,b \in \mathbb{C}^*, a\neq b\,.\]
\indent Now let us apply the Luna \'{E}tale Slice Theorem, to
make a \emph{local} study of $\xi$: in fact it states that if
$Gx$ is a closed semi-stable orbit and $\xi$ is the corresponding
point of $\cq$, then the pointed varieties $(\cq, \xi)$ and
$(N_x//G_x, 0)$ are locally isomorphic in the \'{e}tale topology,
where $N_x=N_{Gx/X,x}$ is the fiber over $x$ of the normal bundle of
$Gx$ in $X$ (for more details about the \'{E}tale Slice Theorem,
see \cite{Luna}, \cite{Walt} and \cite{Drez}).\\
\indent In our case the dimension of the stabilizer
$G_x$ is equal to one and $G_x\cong$ $ \{\textrm{diag}(\lambda^{-2},
\lambda, \lambda), \lambda \in \mathbb{C}^*\} \cong\mathbb{C}^*$.
Moreover the orbit $Gx$ is a 7-dimensional regular variety
in $\mathbb{C}^{12}$ and the space
$T_x\mathbb{C}^{12}=\mathbb{C}^{12}$ can be decomposed
$G_x$-invariantly as the direct sum $T_xGx \oplus N_x$.\\
So we study the action of the torus
$\mathbb{C}^*$ on $N_x$: it is induced by the diagonal
action of $SL_3(\mathbb{C)}$ on $\spa^6(m)$ and it can be written as
\[v_1 \mapsto \lambda^3v_1 ; \quad v_2
\mapsto \lambda^3v_2; \quad v_3 \mapsto \lambda^{-3}v_3; \quad v_4
\mapsto \lambda^{-3}v_4; \quad v_5 \mapsto v_5 \]
where $(v_1,\ldots,v_5)$ is a basis of $N_x\cong \mathbb{C}^5$.\\
In this way a local model of $(\cq, \xi)$ is given by
$( \mathbb{C}^5// \mathbb{C}^*, 0)$ with ``weights''
$(3,3,-3,-3,0)$ that is the 4-dimensional toric variety
\[Y:=\mathbb{C}[T_1,\ldots, T_5]/(T_1T_4-T_2T_3)\,.\]
\indent In conclusion, the variety $(X^{SS}(m)//G, \xi)$, where $\xi$
is a point of the curve $C_{ij}\cong\lin$, is locally isomorphic to the
toric variety $Y$: it is singular and
there are different ways to resolve it (\cite{Ful}, \cite{Cox2}).\\

\subsection{$x_i=x_j$, $x_h=x_l$, $x_k=x_n$}
\indent This study is analogous to the previous one.\\
Consider a polarization $m$ such that it is possible to
``subdivide'' it as $m_i+m_j=m_h+m_l=m_k+m_n$ (for different
indexes); we are examining the configuration $x$, with $x_i=x_j\,,
x_h=x_l\,, x_k=x_n$ (this configuration is a particular
case of the previous one).\\
In the quotient $\cq$ the image of
the orbit $Gx$ is a point $O_{ij,hl,kn}$ that lies on the three
singular curves $C_{ij}, C_{hl}, C_{kn}$.\\
\indent The orbit $Gx$ is minimal, closed and strictly semi-stable:
assume $x$ equal to
\[x=\left(\begin{array}{cccccc}
1& 1& 0& 0& 0& 0\\ 0& 0& 1& 1& 0& 0\\
0& 0& 0& 0&1&1\end{array}\right).\]
Let us apply the \'{E}tale Slice Theorem:
the stabilizer $G_x$ is isomorphic to a 2-dimensional torus
$G_x \cong \{\textrm{diag}(\lambda, \mu, \lambda^{-1}\mu^{-1}),
\lambda,\mu \in \mathbb{C}^*\}$ which implies that $\dim Gx=6$.
By the \'{E}tale Slice Theorem, let us study the action of
$G_x$ on $N_x$: on the basis $\{v_1,\ldots,v_6\}$ of $N_x$ it gives
\[\begin{array}{lll}
v_1 \mapsto \lambda^{-1}\mu\cdot v_1\, ;& \quad v_2 \mapsto \lambda^{-2}\mu^{-1}\cdot v_2\, ;&
\quad v_3 \mapsto \lambda\mu^{-1}\cdot v_3\, ;\\ v_4 \mapsto \lambda^{-1}\mu^{-2}\cdot v_4\, ; &
\quad v_5 \mapsto \lambda^2\mu\cdot v_5\, ; &\quad v_6 \mapsto \lambda\mu^2\cdot v_6\,.
\end{array} \]
\indent It follows that a local model for
$(\cq, O_{ij,hl,kn})$ is given by $Y:=( \mathbb{C}^6//
(\mathbb{C}^*)^2, 0)$, where the action of $(\mathbb{C}^*)^2$ can be
written (in the coordinates $(z_1,\ldots,z_6)$ of $N_x\cong
\mathbb{C}^6$) as
\begin{equation}\label{weights}(\lambda,\mu)(z_1,\ldots,z_6)
\rightarrow(\lambda^{-1}\mu z_1,\lambda^{-2}\mu^{-1}
z_2,\lambda\mu^{-1} z_3,\lambda^{-1}\mu^{-2} z_4,\lambda^2\mu
z_5,\lambda\mu^2 z_6).\end{equation}
Thus we obtain a 4-dimensional toric variety:
\begin{equation}\label{toric quo} Y=\mathbb{C}[T_1,\ldots,
T_5]/(T_1T_2T_3-T_4T_5)\,.\end{equation}
Its singular locus is given
by three lines $s_1=\{(t,0,0,0,0), t \in \mathbb{C}\},$
$s_2=\{(0,t,0,0,0), t \in \mathbb{C}\}$ and $s_3=\{(0,0,t,0,0), t
\in \mathbb{C}\}$ that have a common point, the origin. These lines
correspond to the curves $C_{ij}, C_{hl}, C_{kn}$.\\
\indent A toric representation of $Y$ is determined by a rational,
polyhedral cone $\sigma \subset \mathbb{R}^4$, such that
$Spec(\sigma^\vee \cap \mathbb{Z}^4 )\cong Y$. The generators of the
semi-group $\sigma^\vee \cap \mathbb{Z}^4$ are $w_1,\ldots, w_5 \in
\mathbb{Z}^4$ and satisfy $ w_1+w_2+w_3=w_4+w_5$\,. Assume
\[\begin{array}{c}
w_1=(1,0,0,0)\,,\quad w_2=(0,1,0,0)\,,\quad w_3=(0,0,1,0)\,,\\
w_4=(0,0,0,1)\,,\quad w_5=(1,1,1,-1)\,.\end{array}\] The primitive
elements of $\sigma$ are:
\[ \begin{array}{lll}\mathbf{n}_1=(0,0,1,1)\,,\;& \mathbf{n}_2=(1,0,0,0)\,,\;&
\mathbf{n}_3=(0,0,1,0)\,,\\ \mathbf{n}_4=(0,1,0,1)\,,\;&
\mathbf{n}_5=(1,0,0,1)\,,& \mathbf{n}_6=(0,1,0,0)\,.\end{array}\]
\indent It is clear that the cone $\sigma$ is singular.
\\\indent Let us intesect $\sigma$ with a transversal hyperplane
$\pi$ of $\mathbb{R}^4$ and then consider the projection on $\pi$.
With $\pi: y_1+y_2+y_3+y_4=2$ we get the polytope $\Pi$ of
$\mathbb{R}^3$, with verteces
\[\begin{array}{lll}u_1=(0,0,1)\,,\;& u_2=(2,0,0)\,,\;&
u_3=(0,0,2)\,,\\ u_4=(0,1,0)\,,\;& u_5=(1,0,0)\,,&
u_6=(0,2,0)\,.\end{array}\] \vspace{-1cm}
\begin{figure}[h]
\begin{center}
\includegraphics[bb=15 415 615 770,
scale=0.65]{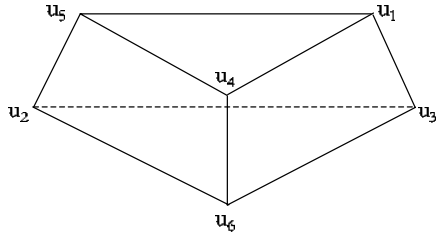} \vspace{-5cm} \caption{Polytope
$\Pi$}\label{fig:Pi}
\end{center}
\end{figure}
\\\indent In conclusion the pointed variety $(\cq, O_{ij,hl,kn})$
is isomorphic to the toric variety $\mathbb{C}[T_1,\ldots,
T_5]/(T_1T_2T_3-T_4T_5)$, where the action has weights
\[\left(\begin{array}{cccccc}-1& -2& 1& -1& 2& 1\\ 1& -1& -1& -2& 1
& 2\end{array}\right)\,.\]

\subsection{$x_h=x_i=x_j$ and $x_k, x_l, x_n$ collinear}
\indent Consider a polarization $m$  such that $m_h+m_i+m_j=|m|/3$
and $m_k+m_l+m_n=2|m|/3$ (for different indexes); then  let us study
the configuration $x$ where: $x_h=x_i=x_j$ and $x_k, x_l, x_n$
collinear.
\\The orbit $Gx$ is minimal, closed,
strictly semi-stable and its image in $\cq$ is a point $O_{hij}$. In
particular $x_k, x_l, x_n$have to be all distinct.\\
As in the previous cases, by the \'{E}tale Slice Theorem, we obtain a local model for $(\cq,
O_{hij})$: this is determined by $Y:=( \mathbb{C}^5// \mathbb{C}^*,
0)$, where the action of $\mathbb{C}^*$ over $\mathbb{C}^5$ with
coordinate $(z_1,\ldots, z_5)$ has weights $(3,\,3,\, 3,\,
3,\,-3)\,$. $Y$ is a 4-dimensional toric variety that corresponds to the smooth affine
variety
\[ Y=\mathbb{C}[T_1,\ldots, T_4]\cong \mathbb{C}^4\,.\]
\indent In conclusion the corresponding point $O_{hij}$ in $\cq$ is nonsingular.\\\\

\indent We have classified the different singularities of $\cq$:
\begin{theo}\label{theo:p26} Let $X=\spa^6$ and $m\in
\mathbb{Z}_{>0}^6$ a polarization:
\begin{enumerate}

\item $m$ s.t.
\begin{itemize}
\item[-] $3 \nmid |m|$,
\item[-] $m_i<|m|/3\,
\forall i$,\end{itemize} then the quotient is geometric;

\item $m$ s.t.
\begin{itemize}
\item[-] $3\mid|m|$,
\item[-] $m_i<|m|/3\, \forall i$,
\item[-] for all couples and triples of indexes we have $m_i+m_j\neq |m|/3$ or $m_h+m_i+m_j\neq
|m|/3$,
\end{itemize}
then the quotient is geometric;

\item $m$ s.t.
\begin{itemize}
\item[-] $3\mid|m|$,
\item[-] there exists an index $i$ s.t. $m_i=|m|/3$,
while for the other indexes $j\neq i,\: m_j<|m|/3\,$,
\end{itemize}
then the quotient is $(\lin)^5(m')//SL_2(\mathbb{C})$; its dimension
is equal to two, and the polarization $m'\in\mathbb{Z}^5_{>0}$ is
obtained from $m$ by eliminating $m_i$;

\item $m$ s.t.
\begin{itemize}
\item[-] $3\mid|m|$,
\item[-] there exist two different indexes $i,j$ s.t. $m_i=m_j=|m|/3$, while for the others
$h\neq i,j,\: m_h<|m|/3\,$,
\end{itemize}
then the quotient is $(\lin)^4(m'')//SL_2(\mathbb{C})\cong \lin $;
the polarization $m''\in\mathbb{Z}^4_{>0}$ is obtained from $m$ by
eliminating $m_i$ and $m_j$;

\item $m$ s.t.
\begin{itemize}
\item[-] $3\mid|m|$,
\item[-]$m_i<|m|/3\, \forall i$,
\item[-] there are two different indexes $i,j$ s.t. $m_i+m_j=|m|/3$,
\end{itemize}
then the quotient is categorical; moreover it includes a curve
$C_{ij}\cong \lin$, that corresponds to strictly semi-stable orbits
s.t. $x_i=x_j$ or $x_h, x_k, x_l, x_n$ collinear. In particular
points $\xi$ of $C_{ij}$ are singular: locally, the variety $(\cq,
\xi)$ is isomorphic to the toric variety\vspace{-0.2cm}
\[\mathbb{C}[T_1,T_2,T_3,T_4,T_5]/(T_1T_4-T_2T_3)\,.\]

\item $m$ s.t.
\begin{itemize}
\item[-] $3\mid|m|$,
\item[-]$m_i<|m|/3\, \forall i$,
\item[-] there is a ``partition'' of $m$ such that $m_i+m_j=m_h+m_l=m_k+m_n$,
\end{itemize}
then the quotient is categorical; moreover it includes three curves
$C_{ij},$ $C_{hl},$ $C_{kn}\cong \lin$, that have a common point
$O_{ij,hl,kn}$. \\In particular $O_{ij,hl,kn}$ is singular: locally
the variety $(\cq, O_{ij,hl,kn})$ is isomorphic to the toric
variety\vspace{-0.2cm}
\[\mathbb{C}[T_1,T_2,T_3,T_4,T_5]/(T_1T_2T_3-T_4T_5)\,.\]

\item $m$ s.t.
\begin{itemize}
\item[-] $3\mid|m|$,
\item[-] $m_i<|m|/3\, \forall i$,
\item[-] there are three indexes $h,i,j$ s.t. $m_h+m_i+m_j=|m|/3$,
\end{itemize}
then the quotient is categorical; moreover it includes a point
$O_{hij}$ that correspond to the minimal, closed, strictly
semi-stable orbit $Gx$ such that $x_h=x_i=x_j$ and
$x_k, x_l, x_n$ are collinear. The point $O_{hij}$ is non singular.\\
\end{enumerate}
\end{theo}

\subsection{Examples}
Now we provide two examples that illustrate how to get explicitly a
quotient, via its coordinates ring, or via an elementary
transformation.

\subsection{$\spa^6(222111)$}
$|m|=9$; by the numerical criterion: $\sum_{k, x_k=y} m_k \leq 3,\,
\sum_{j,x_j \in r} m_j \leq 6$. Then $X^S(m)\subset X^{SS}(m)$\,.\\
Moreover it is easy to verify that there are nine $C_{ij}$ curves,
six $O_{ij,hl,kn}$ points and one $O_{hij}$ point. \\\indent Let us
study the graded algebra of $G$-invariant functions $R_2^6(m)^G$. A
standard tableau $\tau$ of degree $k$ associated to the polarization
$m$ looks like
\begin{equation} \label{eq:tab} \left . \tau=\left[
\begin{array}{ccc} a^1_1 & a^2_2 & a^3_3 \\
a^1_2 & a^2_3 & a^3_4 \\ a^1_3 & a^2_4 & a^3_5 \\ a^1_4 & a^2_5 &
a^3_6 \end{array} \right] \quad \right\}3k \end{equation} where
\[\begin{array}{lll}|a^1_1|=2k,\quad& |a^3_6|=k, \quad&
 |a^1_2|+|a^2_2|=2k, \\|a^1_3|+|a^2_3|+|a^3_3|=2k, &
|a^1_4|+|a^2_4|+|a^3_4|=k, \;&|a^2_5|+|a^3_5|=k,\\  \sum_{i=2}^4
|a^1_i|=k,& \sum_{i=2}^5 |a^2_i|=3k, &\sum_{i=3}^5 |a^3_i|=2k\,.
\end{array}\] Let $ \alpha_3:= |a^1_3|\,, \alpha_4:=|a^1_4|\,,
\beta_3:=|a^3_3|\,, \beta_4:=|a^3_4|$. Then it follows:
\[\begin{array}{lll}
|a^1_1|=2k, & |a^2_2|=k+\alpha_3+\alpha_4, & |a^3_3|=\beta_3,\\
|a^1_2|=k-(\alpha_3+\alpha_4),\; & |a^2_3|=2k-(\alpha_3+\beta_3),\;
&
|a^3_4|=\beta_4, \\
|a^1_3|=\alpha_3,& |a^2_4|=k-(\alpha_4+\beta_4),&
|a^3_5|=2k-(\beta_3+\beta_4),\\
|a^1_4|=\alpha_4,& |a^2_5|=\beta_3+\beta_4-k,&|a^3_6|=k.
\end{array}\]
Moreover $\alpha_3, \alpha_4, \beta_3, \beta_4$ must satisfy the
inequalities:
\[
\begin{array}{lll} 0\leq \alpha_3,\alpha_4, \beta_3,\beta_4 \leq 2k,&
\alpha_3+2\alpha_4 \leq \beta_3, & \alpha_3+\alpha_4 \leq k, \\
k+\alpha_4\leq \beta_3+\beta_4 \leq 2k,\;& \beta_3 \leq
k+\alpha_3+\alpha_4,\; & 2\beta_3+\beta_4 \leq 3k+\alpha_4.
\end{array}\] Assume \[ x:=\alpha_4,\quad
y:=\alpha_3+\alpha_4,\quad z:=\beta_3,\quad w:=\beta_3+\beta_4\,;\]
the standard tableau $\tau$ (\ref{eq:tab}) is completely determined
by the vector $(x,y,z,w)$ that satisfy:
\[\begin{array}{cccc} 0 \leq x \leq y \leq k,\;& 0 \leq z \leq w \leq 2k,
& 0\leq y+z-x \leq 2k,& \\
x+y \leq z \leq y+k,& z\leq w \leq k+z,& 0\leq w+x-z \leq k,\; & w
\geq x+k\,.
\end{array}\]
After few calculations we find out that for any $k$, there are
\[\frac{1}{8}(k^4+6k^3+15k^2+18k)+1\,(=\,\dim(R^6_2(m)_k^G)\,)\]
standard tableaux. Thus the Hilbert function of the graded ring
$R^6_2(m)^G$ is equal to
\[\sum_{k=0}^\infty
\left(\frac{1}{8}(k^4+6k^3+15k^2+18k)+1\right)t^k =
\frac{1-t^3}{(1-t)^6}\,.\]  This suggests that the quotient $\cq$ is
isomorphic to a cubic hypersurface in $\mathbb{P}^5(\mathbb{C})$.
\\ \indent First of all we have the following generators of $R^6_2(m)^G$\,:
\[\begin{array}{lll}t_0=[124][135][236], \;& t_1=[123][135][246],
\;& t_2=[123][134][256],\\ t_3=[123][125][346],\;&
t_4=[123][124][356], \;& t_5=[123][123][456].\end{array}\] For every
$(i,j)\neq (2,3),(3,2)$, the product $t_it_j$ is a standard tableau
function from $R^6_2(m)_2^G$. Applying the straightening algorithm
(that allows to write any tableau function as a linear combination
of tableau standard functions), we obtain:
\begin{equation}\label{eq:u}
t_2t_3=t_1t_4-u+t_5(-t_0+t_1-t_2-t_3+t_4-t_5).\end{equation} So the
standard monomial $u=[123][123][123][145][246][356]$ can be
expressed as polynomials of degree two in the $t_i$.
\\\indent In we take a tableau function $\mu_{(x,y,z,w,k)}$ corresponding to
a standard tableau $\tau$ (\ref{eq:tab}), we can write it as
\[\mu_{(x,y,z,w,k)} = \left\{\begin{array}{ll}
t_0^{k+x-z}t_1^{k+z-x-w}t_2^{w-y-k}t_4^{y-x}t_5^x,  & z\leq x+k, \;w \leq k+z-x;\\
t_0^{k+x-z}t_1^{z-x-y}t_3^{k+y-w}t_4^{w-x-k}t_5^x, & z\leq x+k,\; y\leq z-x;\\
t_1^{3k+x-w-z}t_2^{w-y-k}t_4^{k+y-z}t_5^{x}u^{z-x-k},\: & z\geq x+k,\; w\leq 3k+x-z;\\
t_1^{2k+x-y-z}t_3^{k+y-w}t_4^{w-z}t_5^{x}u^{z-x-k}, & z\geq x+k,\;
y\leq 2k+x-z.
\end{array} \right.
\]
\indent Applying the straightening algorithm to the non-standard
product $t_0u$, we have:
\[t_0u\,=\,t_1t_4(t_1-t_2-t_3+t_4-t_5)\,.\]
Then by relation (\ref{eq:u}), it follows
\[t_0\big(t_1t_4-t_2t_3+t_5(-t_0+t_1-t_2-t_3+t_4-t_5)\big)=
t_1t_4(t_1-t_2-t_3+t_4-t_5)\, \Rightarrow\]
\[t_0\big(-t_2t_3+t_5(-t_0+t_1-t_2-t_3+t_4-t_5)\big)=
t_1t_4(-t_0+t_1-t_2-t_3+t_4-t_5)\, \Rightarrow\]
\[(-t_0+t_1-t_2-t_3+t_4-t_5)(t_0t_5-t_1t_4)-t_0t_2t_3=0\]
Let \begin{equation}\label{eq:cub} F_3=
(-T_0+T_1-T_2-T_3+T_4-T_5)(T_0T_5-T_1T_4)-T_0T_2T_3\,,\end{equation}
there is a surjective homomorphism of the graded algebras
\[\mathbb{C}[T_0,T_1,T_2,T_3,T_4,T_5]/(F_3(T_0,T_1,T_2,T_3,T_4,T_5))
\longrightarrow R^6_2(m)^G\,.\] \indent Thus the quotient $\cq$ is
isomorphic to the cubic hypersurface $F_3(T_0,T_1,T_2,T_3,T_4,T_5)=0$.

\subsection{$\spa^6(221111)$}
$|\mw|=8$; by the numerical criterion $\sum_{k, x_k=y} \mw_k \leq
8/3$, $\sum_{j,x_j \in r} \mw_j \leq 16/3$ and thus $X^S(\mw)=
X^{SS}(\mw)$\,.
\\\indent
In order to determine this geometric quotient, we have to introduce
the elementary transformation
$\mw=(221111)\stackrel{+1_3}{\longrightarrow}(222111)=m$, and
consequently \[\widehat{\theta}: X^S(\mw)/G\longrightarrow \cq\,.\]
\indent First of all let us study
$\widehat{\theta}^{-1}(O_{456})$: by relation (\ref{dim fiber i in
J'}) its dimension is equal to $d=3$; the
semi-stable orbits of $\ssm$ that determine $O_{456}$ in the
quotient $\cq$ and are included in $X^S(\mw)$, are characterized by
$x_1, x_2, x_3$ collinear. Applying a projectivity of $\spa$ such
that it fixes the line that contains $x_1, x_2, x_3$\,, we have
$\widehat{\theta}^{-1}(O_{456})\cong\mathbb{P}^3(\mathbb{C})$.
\\\indent Then $\widehat{\theta}^{-1}(\xi), \xi \in C_{ij}$; studying
how semi-stable orbits change going from $\ssm$ to $X^S(\mw)$, there
can be two different cases: coincidence or collinearity.
\begin{enumerate}
\item Consider the curve $C_{14}$: by the numerical
criterion for $X^S(\mw)$, orbits which have $x_2, x_3, x_5, x_6$
collinear are stable. In particular by relation (\ref{dim fiber i in
J'}), the dimension of $\widehat{\theta}^{-1}(\xi_1)$, $\xi_1 \in
C_{14}$ is equal to $d=1$: in fact
\begin{equation}\label{resol C14}\widehat{\theta}^{-1}(\xi_1) \cong
\lin\,.\end{equation}
\item Consider the curve $C_{36}$:
by the numerical criterion for $X^S(\mw)$ orbits which have
$x_3=x_6$ are stable. In particular by relation (\ref{dim fiber i in
J}), the dimension of $\widehat{\theta}^{-1}(\xi_2)$, $\xi_2 \in
C_{36}$ is equal to $d=1$; in fact
\begin{equation}\label{resol C36}\widehat{\theta}^{-1}(\xi_2) \cong
\lin\,.\end{equation}
\end{enumerate}

\begin{figure}[h]
\begin{center}
\includegraphics[bb=0 540 650 770,
scale=0.5]{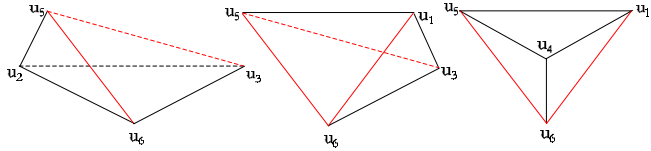} \vspace{-2cm} \caption{The three polytopes of
the subdivision (221111) of $\Pi$}\label{fig: 3Pi}
\end{center}
\end{figure}

\indent Let us study $\widehat{\theta}^{-1}(O_{ij,hl,kn})$; consider
$O_{14,25,36}$. Strictly semi-stable orbits that contain the orbit
$Gx$ $(x_1=x_4, x_2=x_5, x_3=x_6)$ in their closure, are
characterized by one of the following properties:\\
\indent $1.$ $x_1=x_4$ and $x_1, x_2, x_5$ collinear;\indent $2.$ $x_1=x_4$ and $x_1, x_3, x_6$ collinear;\\
\indent$3.$ $x_2=x_5$ and $x_1, x_2, x_4$ collinear;\indent $4.$ $x_2=x_5$ and $x_2, x_3, x_6$ collinear;\\
\indent$5.$ $x_3=x_6$ and $x_1, x_3, x_4$ collinear;\indent $6.$ $x_3=x_6$ and $x_2, x_3, x_5$ collinear.\\
In particular configurations $1,2,3,4$ are unstable for the
polarization $\mw$, while $5$ and $6$ are included in $X^S(\mw)$;
moreover these sets have a common configuration: $(x_3=x_6$, $x_1,
x_3, x_4$ collinear, $x_2, x_3, x_5$ collinear$)$:
\setlength{\unitlength}{1cm}
\[\begin{picture}(3,3)
\put(0.8,3){\line(1,-3){1}} \put(1.2,3){\line(-1,-3){1}}
\put(1,2.4){\circle*{0.15}} \put(1.3,2.3){$x_3=x_6$}
\put(0.7,1.5){\circle{0.15}} \put(0.2,1.4){$x_1$}
\put(1.3,1.5){\circle{0.15}} \put(1.5,1.4){$x_2$}
\put(0.45,0.8){\circle{0.15}} \put(0,0.7){$x_4$}
\put(1.55,0.8){\circle{0.15}} \put(1.75,0.7){$x_5$}
\end{picture}\] Every one of these two sets of stable configurations
determine a copy of $\lin$ in the quotient $X^S(\mw)/G$: thus these
two copies of $\lin$ have a common point.
\[\widehat{\theta}^{-1}(O_{ij,hl,kn})\cong \lin\cup \lin\;
\textrm{with a common point}\,.\]
\indent We can get this result in
a different way, by constructing a subdivision of the polytope $\Pi$
(figure \ref{fig:Pi}).
\\\indent Since $X^{US}(m) \subset X^{US}(\mw)$ and $\left(
X^{US}(\mw) \setminus X^{US}(m) \right) \subset X^{SSS}(m)$, we
determine (locally in $N_x$), which strictly
semi-stable orbits for the polarization $m$ are unstable for $\mw$.
By the machinery of the theory of homogeneous coordinates for a
toric variety (\cite{Cox1},\cite{Cox2}, \cite{Dol}),
the local resolution of $(\cq, O_{14,25,36}) \cong
(\mathbb{C}^6/(\mathbb{C}^*)^2, 0)$ in the quotient $X^S(\mw)/G$
is determined by $ (\mathbb{C}^6 \setminus Z)//H\,,$ where
$\mathbb{C}^6 \setminus  Z= \mathbb{C}^6\setminus \{z \in
\mathbb{C}^6\,| z_1z_4=0, z_2z_3=0, z_2z_4=0\}$\,, and $H$ is the 2-dimensional torus
$H=\{(\lambda_1,\lambda_2,\lambda_1^{-1}, \lambda_1^{-1}\lambda_2,
\lambda_2^{-1}, \lambda_1\lambda_2^{-1}), \lambda_1,\lambda_2 \in
\mathbb{C}^* \}$.
\\\indent The set $\mathbb{C}^6
\setminus Z$ describes a particular resolution of $\Pi$.
\begin{figure}[h]
\begin{flushright}
\includegraphics[bb=25 415 590 750,
scale=0.65]{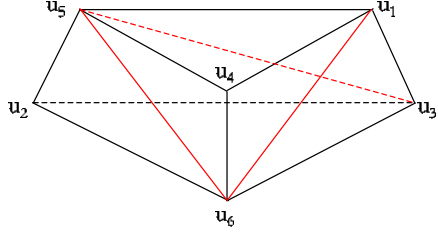}\vspace{-4.5cm} \caption{Subdivision of type
(221111) of $\Pi$}\label{fig: SudPi}
\end{flushright}
\end{figure}
\\We can find three simplicial polytopes: figure \ref{fig: 3Pi}. \\
\\
\indent The toric representation of $Y$, described by the polytope
$\Pi$, is determined by the cone $\sigma$: to solve its
singularities let us construct a fan $\Sigma$, refinement of
$\sigma$. By the theory of toric varieties, there exists a proper,
birational morphism $\varphi$
\[X_{\Sigma}\cong (\mathbb{C}^6 \setminus Z)//H\cong (\mathbb{C}^6 \setminus Z)//
(\mathbb{C}^*)^2 \stackrel{\varphi}{\longrightarrow} (\mathbb{C}^6
//(\mathbb{C^*)}^2)\cong(N_x//G_x)\cong X_\sigma,\] induced by
the identity over the lattice $\mathbb{R}^4$: this application
allows us to specify the map $\widehat{\theta}$:
\[\widehat{\theta}: X^{S}(\mw)/G \longrightarrow
X^{SS}(m)//G\,.\]
\indent First of all let us take a cover of $(\mathbb{C}^6 \setminus
Z)$: for example the three open sets $U_1, U_2, U_3$:
\[\begin{array}{c}U_1=\mathbb{C}^6 \setminus \{z \in \mathbb{C}^6\,|\;
z_1z_4=0\};\qquad
U_2=\mathbb{C}^6 \setminus \{z \in \mathbb{C}^6\,|\; z_2z_3=0\};\\
U_3=\mathbb{C}^6 \setminus \{z \in \mathbb{C}^6\,|\; z_2z_4=0\}.
\end{array}\] Now let us consider the action of $H\cong(\mathbb{C^*})^2$
on these three open sets and construct the three quotients: in
the first case, the quotient $\widetilde{U_1}=U_1//H$
is the smooth variety $\mathbb{C}[X_1,X_2,X_3,X_4,X_6]/(X_2-X_4X_6)$.\\
In the same way
$\widetilde{U}_2=U_2//H=\mathbb{C}[Y_1,Y_2,Y_3,Y_5,Y_7]/(Y_3-Y_5Y_7)$
and
$\widetilde{U}_3=U_3//H=\mathbb{C}[Z_1,Z_2,Z_3,Z_8,
Z_9]/(Z_1-Z_8Z_9)$.\\\\
How do these quotients $\widetilde{U}_i (i=1,2,3)$ fit together? We
have the following ``gluing''
\begin{equation}\label{eq:glue}\begin{array}{lll}X_1=Y_1=Z_8Z_9 &
Y_1=X_1=Z_8Z_9 & Z_2=X_4X_6=Y_2\\
X_3=Y_5Y_7=Z_3 & Y_2=X_4X_6=Z_2& Z_3=X_3= Y_5Y_7\\
X_4=Y_1Y_2Y_7=Z_2Z_8 & Y_5=X_1X_3X_6=Z_3Z_8& Z_8=X_6^{-1}=Y_1Y_7\\
X_6=(Y_1Y_7)^{-1}=Z_8^{-1}& Y_7=(X_1X_6)^{-1}=Z_9^{-1}&
Z_9=X_1X_6=Y_7^{-1}
\end{array}\end{equation} \indent The birational maps
$\widehat{\theta}_i: \widetilde{U}_i \rightarrow Y$ that resolve the
singularities of $Y$ are described by the pull back of the generators of the ring of
$G_x$-invariant functions $(T_1, T_2, T_3,T_4,T_5)$:
\[\begin{array}{lll}\widehat{\theta}_1^*(T_1)=X_1,& \widehat{\theta}_2^*(T_1)=Y_1,&
\widehat{\theta}_3^*(T_1)=Z_8Z_9,\\
\widehat{\theta}_1^*(T_2)=X_4X_6,& \widehat{\theta}_2^*(T_2)=Y_2,&
\widehat{\theta}_3^*(T_2)=Z_2,
\\ \widehat{\theta}_1^*(T_3)=X_3,&  \widehat{\theta}_2^*(T_3)=Y_5Y_7,&
\widehat{\theta}_3^*(T_3)=Z_3,\\
\widehat{\theta}_1^*(T_4)=X_4, &
\widehat{\theta}_2^*(T_4)=Y_1Y_2Y_7,&
\widehat{\theta}_3^*(T_4)=Z_2Z_8,\\
\widehat{\theta}_1^*(T_5)=X_1X_3X_6, &
\widehat{\theta}_2^*(T_5)=Y_5,& \widehat{\theta}_3^*(T_5)=Z_3Z_9.
\end{array}\] The point $O_{14,25,36}$ corresponds to the origin
in $Y$: let us study $\widehat{\theta}_i^{-1}(0)$
\[ \begin{array}{c}\widehat{\theta}_1^{-1}(0)=(0,0,0,t_1)\cong
\mathbb{C},\quad \widehat{\theta}_2^{-1}(0)=(0,0,0,u_1)\cong
\mathbb{C},
\\ \widehat{\theta}_3^{-1}(0)=(0,0,t_2,u_2) \cong \mathbb{C} \cup \mathbb{C}
\end{array}\]
where  $t_1,u_1,t_2,u_2 \in \mathbb{C}$ and $t_2u_2=0$. \\In
particular the fiber $\widehat{\theta}_3^{-1}(0)$ is isomorphic to
the union of two copies of $\mathbb{C}$ that have a common point
$(0,0,0,0) \in \widetilde{U}_3$. Moreover by the gluing (\ref{eq:glue}), $t_1,
t_2 \in \mathbb{C}$ give a cover of $\lin$, just like $u_1, u_2 \in
\mathbb{C}$.\\\indent In conclusion the resolution of
$O_{14,25,36}$ in $X^S(221111)/G$ is determined by the union of two
copies of $\lin$ that have a common point
\[\widehat{\theta}^{-1}(O_{14,25,36})\cong \lin\cup
\lin\,\quad \textrm{with a common point.}\]
\indent Let us calculate the resolutions of the three singular
curves $C_{14}, C_{25}, C_{36}$ that meet in $O_{14,25,36}$: we know
that there is a correspondence between $C_{ij}, C_{hl}, C_{kn}$ and
the three lines $s_3=\{(0,0,t,0,0)\}, s_2=\{(0,t,0,0,0)\},
s_1=\{(t,0,0,0,0)\}$ of $Y$. Now let us calculate the
fiber of a ``generic'' point of each line $s_j$, for the maps
$\widehat{\theta}_i$.\\
Let $\xi_3 \in C_{14}$:
$\widehat{\theta}_1^{-1}(\xi_3)=(0,t,0,\tau)$,
$\widehat{\theta}_2^{-1}(\xi_3)=\textrm{Imposs.}$,
$\widehat{\theta}_3^{-1}(\xi_3)=(0,t,\tau^{-1},0)$; thus
\[\widehat{\theta}^{-1}(\xi_3)\cong \lin\,,\quad
\forall \xi_3\in C_{14}\quad \xi_3\neq O_{ij,hl,kn}.\]\\
In the same way for $\xi_2
\in C_{25}$ and $\xi_1 \in C_{36}$, $\xi_1,\xi_2 \neq O_{ij,hl,kn}$ we obtain:
\[\widehat{\theta}^{-1}(\xi_2)\cong \lin\,,\quad
\widehat{\theta}^{-1}(\xi_1)\cong \lin\,.\]
\indent In conclusion the map \[\widehat{\theta}:
X^S(\mw)/G=(\mathbb{P}^2)^6(221111)/G \longrightarrow
(\mathbb{P}^2)^6(222111)//G=\cq\] determines the quotient
$X^S(\mw)/G$: in fact $\widehat{\theta}$ is an isomorphism over
\[X^S(\mw)/G \setminus \left( \bigcup_{\xi \in
S}\widehat{\theta}^{-1}(\xi) \right)\stackrel{\sim}{\longrightarrow}
\gq\,,\] where $S=\{\xi \in X^{SSS}(m)//G\}$.\\ Then the map
$\widehat{\theta}$ is a contraction of subvarieties over
$\bigcup_{\xi \in S}\widehat{\theta}^{-1}(\xi)$:
\begin{itemize}
\item[-] if $\xi \in C_{ij}$,
then $\widehat{\theta}^{-1}(\xi)=\lin$;
\item[-] if $\xi=O_{ij,hl,kn}$, then $\widehat{\theta}^{-1}(\xi)=
\lin\cup \lin$, with a common point;
\item[-] if $\xi=O_{456}$, then $\widehat{\theta}^{-1}(\xi)=
\mathbb{P}^3(\mathbb{C})$.
\end{itemize}


\noindent \textsc{Francesca Incensi}\\
Dipartimento di Matematica, Universit\`{a} di Bologna, Italy\\
E-mail address: \textsl{incensi@dm.unibo.it}
\end{document}